\documentclass[letterpaper, 10 pt, conference]{ieeeconf} 
\IEEEoverridecommandlockouts                              
\usepackage{graphicx}
\usepackage{amsmath}
\usepackage{dsfont}
\usepackage{xcolor}
\usepackage{amsfonts}
\usepackage{amssymb}
\usepackage{epsfig,psfrag,pstricks}
\usepackage{subfigure}
\usepackage{tikz}
\usepackage{tkz-berge}
\usetikzlibrary{shapes,arrows}
\DeclareMathOperator{\im}{im}

\DeclareMathOperator{\spa}{span}
\DeclareMathOperator{\sat}{sat}

\newtheorem{theorem}{Theorem}
\newtheorem{lemma}[theorem]{Lemma}

\newtheorem{corollary}[theorem]{Corollary}
\newtheorem{exmp}{Example}[section]
\newtheorem{remark}[theorem]{Remark}

\title{\LARGE \bf
Stability of dynamical distribution
networks with arbitrary flow constraints and unknown in/outflows*
}

\author{Jieqiang Wei$^{1}$ and Arjan van der Schaft$^{2}$
\thanks{*The work of the first author is supported by the Chinese Science
Council (CSC). The research of the second author leading to these
results has received funding from the EU 7th Framework Programme
[FP7/2007-2013] under grant agreement no. 257462 HYCON2 Network of
Excellence.}
\thanks{$^{1}$,$^{2}$Johann Bernoulli Institute for Mathematics and Computer
Science, University of Groningen, PO Box 407, 9700 AK, the
Netherlands. {\tt\small J. Wei@rug.nl, A.J.van.der.Schaft@rug.nl}}%
}

\begin{document}

\maketitle
\thispagestyle{empty}
\pagestyle{empty}

\begin{abstract}

A basic model of a dynamical distribution network is considered, modeled as a
directed
graph with storage variables corresponding to every vertex and flow inputs
corresponding to every edge, subject to unknown but constant
inflows and outflows. We analyze the dynamics of the system in closed-loop with a distributed proportional-integral controller structure, where the
flow inputs are constrained to take value in closed intervals. Results from our previous work are extended to general flow constraint intervals, and conditions for asymptotic load balancing are derived that rely on the
structure of the graph and its flow constraints.

\end{abstract}

\section{INTRODUCTION}

In this paper we study a basic model for the dynamics of a distribution
network. Identifying the network with a directed graph we associate with every
vertex of the graph a state variable corresponding to {\it storage}, and with
every edge a control input variable corresponding to {\it flow}, which is
constrained
to take value in a given closed interval. Furthermore, some of the vertices
serve as terminals where an unknown but constant flow may enter or leave the
network in such a way that the total sum of inflows and outflows is equal to
zero. The control problem to be studied is to derive necessary and sufficient
conditions for a distributed control structure (the control input corresponding
to a
given edge only depending on the difference of the state variables of the
adjacent vertices) which will ensure that the state variables associated to all
the vertices will converge to the same value equal to the average of the initial
condition,
irrespective of the values of the constant unknown inflows and outflows.

The structure of the paper is as follows. Preliminaries and notations will
be given in Section 2. In Section 3 we will briefly recall how in the absence of
constraints on the flow input variables a distributed proportional-integral (PI)
controller structure,
associating with every edge of the graph a controller state, will solve the
problem if and only if the graph is weakly connected; see also \cite{Wei2012}.
This will
be shown by identifying the closed-loop system as a port-Hamiltonian system,
with state variables associated both to the vertices and the edges of the graph,
in line with the general definition of port-Hamiltonian systems on graphs
\cite{schaftSIAM, schaftCDC08, schaftBosgrabook,
schaftNECSYS10}; see also \cite{allgower11,Mesbahi11}.

In Sections 4 and 5 the same problem is studied in the presence of constraints
on the flow inputs. In \cite{Blanchini2000}, the authors consider a similar
model and present a discontinuous Lyapunov-based controller to stabilize the
system without violating the storage and flow constraints. In \cite{Bauso2011}, 
using the same model as in \cite{Blanchini2000}, the authors focus on a
different problem of driving the state to a small neighborhood of the reference
value and relate the control input value at equilibrium to an optimization problem. In the current paper we will
generalize most of the results of our previous work \cite{Wei2013} to
the case of {\it arbitrary} constraint intervals, making use of a new technique
extending the graph to a graph with a larger number of edges admitting a
coverage by non-overlapping cycles. Section 6 contains the conclusions.

\section{Preliminaries and notations}

First we recall some standard definitions regarding directed graphs,
as can be found e.g. in \cite{Bollobas98}. A \textit{directed graph}
$\mathcal{G}$ consists of a finite set $\mathcal{V}$ of \textit{vertices}
and a finite set $\mathcal{E}$ of \textit{edges}, together
with a mapping from $\mathcal{E}$ to the set of ordered pairs of
$\mathcal{V}$, where no self-loops are allowed. Thus to any edge
$e\in\mathcal{E}$ there corresponds an ordered pair
$(v,w)\in\mathcal{V}\times\mathcal{V}$
(with $v\not=w$), representing the tail vertex $v$ and the head
vertex $w$ of this edge.

A directed graph is specified by its \textit{incidence
matrix} $B$, which is an $n\times m$ matrix, $n$ being the
number of vertices and $m$ being the number of edges, with $(i,j)^{\text{th}}$
element equal to $1$ if the $j^{\text{th}}$ edge is towards vertex
$i$, and equal to $-1$ if the $j^{\text{th}}$ edge is originating from
vertex $i$, and $0$ otherwise. Since we will only consider directed
graphs in this paper `graph' will throughout mean `directed graph'
in the sequel.
A directed graph is {\it strongly connected} if it is
possible to reach any vertex starting from any other vertex by traversing edges
following their directions. A directed graph is called {\it weakly connected}
if it is possible to reach any vertex from every other vertex using the edges
{\it not} taking into account their direction. A graph is weakly connected if
and only if $\ker B^T = \spa \mathds{1}_n$. Here $\mathds{1}_n$ denotes the
$n$-dimensional vector with all elements equal to $1$. A graph that is not
weakly connected falls apart into a number of weakly connected subgraphs, called
the connected components. The number of connected components is equal to $\dim
\ker B^T$.
For each vertex, the number of incoming edges is called the {\it in-degree} of
the vertex and the number of outgoing edges its out-degree. A graph is called
{\it balanced} if and only if the in-degree and out-degree of every vertex are
equal. A graph is balanced if and only if $\mathds{1}_n \in \ker B$.

Given a graph, we define its \textit{vertex space} as the vector space of all
functions from $\mathcal{V}$ to some linear space $\mathcal{R}$. In the rest of
this paper we will take $\mathcal{R}=\mathbb{R}$, in which case
the vertex space can be identified with $\mathbb{R}^{n}$. Similarly, we define
its \textit{edge space} as the
vector space of all functions from $\mathcal{E}$ to $\mathcal{R} = \mathbb{R}$,
which can be identified with $\mathbb{R}^{m}$. In this way, the incidence matrix
$B$ of the graph can be also regarded as the matrix representation of a linear
map from the edge space $\mathbb{R}^m$ to the vertex space $\mathbb{R}^n$.

\noindent
{\bf Notation}: For $a,b\in\mathbb{R}^m$ the notation $a \leqslant b$ (resp. $ A < b$) will
denote elementwise inequality $a_i \leq b_i$ (resp. $a_i < b_i$),\,$i=1,\ldots,m$. For $a <
b$ the multidimensional
saturation function
$\sat(x\,;a,b): \mathbb{R}^m\rightarrow\mathbb{R}^m$ is defined as
\begin{equation}
\sat(x\,;a,b)_i  = \left\{ \begin{array}{ll}
a_i & \textrm{if $x_i\leq a_i,$}\\
x_i & \textrm{if $a_i<x_i<b_i,$}\\
b_i & \textrm{if $x_i\geq b_i$},
\end{array}
\, i=1,\ldots,m. \right.
\end{equation}

\section{A dynamical network model with PI controller}
Consider the following dynamical system defined on the graph
$\mathcal{G}$
\begin{equation}\label{system}
\begin{array}{rcl}
\dot{x} & = & Bu, \quad x \in \mathbb{R}^n, u \in \mathbb{R}^m \\[2mm]
y & = & B^T \frac{\partial H}{\partial x}(x), \quad y \in \mathbb{R}^m,
\end{array}
\end{equation}
where $H: \mathbb{R}^n \to \mathbb{R}$ is a differentiable function, and
$\frac{\partial H}{\partial
x}(x)$ denotes the column vector of partial derivatives of $H$. Here the
$i$-th element $x_i$ of the state vector $x$ is the state variable
associated to the $i$-th vertex, while $u_j$ is a flow input variable
associated to the
$j$-th edge of the graph. System (\ref{system}) defines a
{\it port-Hamiltonian system} (\cite{vanderschaftmaschkearchive,
vanderschaftbook}), satisfying the energy-balance
\begin{equation}
\frac{d}{dt}H = u^Ty.
\end{equation}
Note that geometrically its state space is the vertex space, its input space is
the edge space, while its output space is the dual of the
edge space \cite{schaftSIAM}. Many distribution networks are of this form; see 
\cite{Wei2012}, \cite{schaftSIAM} for further background.

Furthermore, we extend the dynamical system (\ref{system}) with a vector
$d$ of {\it inflows and outflows}
\begin{equation}\label{system1}
\begin{array}{rcl}
\dot{x} & = & Bu + Ed, \quad x \in \mathbb{R}^n, u \in \mathbb{R}^m, \quad d \in
\mathbb{R}^k \\[2mm]
y & = & B^T \frac{\partial H}{\partial x}(x), \quad y \in \mathbb{R}^m,
\end{array}
\end{equation}
where $E$ is an $n \times k$ matrix whose columns consist of exactly one entry equal
to $1$ (inflow) or $-1$ (outflow), while the rest of the elements is zero. Thus
$E$ specifies the $k$ terminal vertices where flows can enter or leave the
network ({\it sources} and {\it sinks}).

As in \cite{Wei2012}, \cite{Wei2013} we will regard $d$ as a vector of {\it
constant disturbances},
and we want to investigate control schemes which
ensure asymptotic load balancing of the state vector $x$ irrespective of the
unknown value of $d$. The
simplest strategy is to apply a proportional output feedback (as in \cite{Bauso2011})
\begin{equation}\label{Pcontroller}
u = - Ry = -RB^T\frac{\partial H}{\partial x}(x),
\end{equation}
where $R$ is a diagonal matrix with strictly positive diagonal elements
$r_1,\ldots,r_m$. Note that this defines a {\it decentralized} control scheme if
$H$ is of the form $H(x)= H_1(x_1) + \ldots + H_n(x_n)$, in which case the
$i^{\text{th}}$ input is given as $r_i$ times the difference of the component of
$\frac{\partial H}{\partial x}(x)$ corresponding to the head vertex of the
$i^{\text{th}}$ edge and the component of $\frac{\partial H}{\partial x}(x)$
corresponding to its tail vertex.

However, for $d \neq 0$ proportional control $u=
-Ry$ will not be sufficient to reach load balancing, since the disturbance $d$
can only be attenuated at the
expense of increasing the gains in the matrix $R$. Hence we
consider instead a {\it proportional-integral} (PI) control structure, given by
\begin{equation}\label{PI}
\begin{array}{rcl}
\dot{x}_c & = & y ,\\[2mm]
u & = &-Ry - \frac{\partial H_c}{\partial x_c}(x_c),
\end{array}
\end{equation}
where $H_c(x_c)$ denotes the storage function (energy) of the controller. Note
that this PI controller is of the same distributed nature as the static output
feedback $u=-Ry$.

The $j$-th element of the controller state $x_c$ can be regarded as an
additional state variable corresponding to the $j$-th edge. Thus $x_c
\in
\mathbb{R}^m$, the edge space of the network. The closed-loop system resulting
from the PI control (\ref{PI}) is given as
\begin{equation}\label{closedloop}
\begin{bmatrix} \dot{x} \\[2mm] \dot{x}_c \end{bmatrix} =
\begin{bmatrix} -BRB^T & -B \\[2mm] B^T & 0 \end{bmatrix}
\begin{bmatrix} \frac{\partial H}{\partial x}(x) \\[2mm] \frac{\partial
H_c}{\partial x_c}(x_c) \end{bmatrix} +
\begin{bmatrix} E \\[2mm] 0 \end{bmatrix} d,
\end{equation}
This is again a port-Hamiltonian system,
with total
Hamiltonian $$H_{\mathrm{tot}}(x,x_c)\\ := H(x) + H_c(x_c),$$ and satisfying
the
energy-balance
\begin{equation}\label{Lyapunov}
\frac{d}{dt}  H_{\mathrm{tot}}= - \frac{\partial^T H}{\partial x}(x)BRB^T
\frac{\partial H}{\partial x}(x) + \frac{\partial^T H}{\partial x}(x)Ed
\end{equation}
Consider now a constant disturbance $\bar{d}$ for which there exists a {\it
matching}
controller state $\bar{x}_c$, i.e.,
\begin{equation}\label{matching}
E \bar{d} = B\frac{\partial H_c}{\partial x_c}(\bar{x}_c).
\end{equation}
By modifying the total Hamiltonian $H_{\mathrm{tot}}(x,x_c)$
into the candidate Lyapunov function
\begin{equation}
\begin{aligned}
V_{\bar{d}}(x,x_c) :=& H(x) + H_c(x_c) \\ &- \frac{\partial^T H_c}{\partial
x_c}(\bar{x}_c)(x_c - \bar{x}_c) - H_c(\bar{x}_c),
\end{aligned}
\end{equation}
the following theorem is obtained \cite{Wei2012, Wei2013}.
\begin{theorem}
Consider the system (\ref{system1}) on the graph $\mathcal{G}$ in closed-loop
with the PI-controller (\ref{PI}).
Let the constant disturbance $\bar{d}$ be such that there exists a $\bar{x}_c$
satisfying the matching equation (\ref{matching}). Assume that
$V_{\bar{d}}(x,x_c)$ is radially unbounded. Then the trajectories of the
closed-loop system (\ref{closedloop}) will converge to an element of the load
balancing set
\begin{equation}
 \mathcal{E}_{\mathrm{tot}} = \{ (x,x_c) \mid \frac{\partial H}{\partial x}(x)=
 \alpha \mathds{1}, \, \alpha \in \mathbb{R}, \, B\frac{\partial H_c}{\partial
 x_c}(x_c) = E\bar{d}\, \}.
\end{equation}
if and only if $\mathcal{G}$ is weakly connected.
\end{theorem}
\begin{corollary}
If $\ker B = 0$, which is equivalent (\cite{Bollobas98}) to the graph having no
{\it cycles}, then for every $\bar{d}$ there exists a {\it unique} $\bar{x}_c$
satisfying (\ref{matching}), and convergence is towards the set
$\mathcal{E}_{\mathrm{tot}} = \{ (x, \bar{x}_c)
\mid \frac{\partial H}{\partial x}(x) = \alpha \mathds{1}, \alpha \in
\mathbb{R}, \, x_c=\bar{x}_c \}$.
\end{corollary}
\begin{corollary}
In case of the standard quadratic Hamiltonians $H(x) = \frac{1}{2} \| x \|^2,
H_c(x_c)=\frac{1}{2} \| x_c \|^2$ there exists for every $\bar{d}$
a controller state $\bar{x}_c$ such that (\ref{matching}) holds if and only if
\begin{equation}\label{matching1}
\im E \subset \im B.
\end{equation}
Furthermore, in this case $V_{\bar{d}}$ equals the radially unbounded function
$\frac{1}{2}
\| x \|^2 + \frac{1}{2} \| x_c - \bar{x}_c \|^2$, while convergence will be
towards the load balancing set $\mathcal{E}_{\mathrm{tot}} = \{ (x,x_c) \mid x =
\alpha
\mathds{1}, \alpha \in \mathbb{R},\, Bx_c = E\bar{d}\}$.
\end{corollary}

A necessary (and in case the graph is weakly connected necessary {\it and}
sufficient) condition for the inclusion $\im E \subset \im B $ is that
$\mathds{1}^TE =
0$. In its turn $\mathds{1}^TE =
0$ is equivalent to the fact that for every $\bar{d}$ the total inflow into the
network equals to the
total outflow). The condition $\mathds{1}^TE = 0$ also implies
\begin{equation}
\mathds{1}^T\dot{x} = -\mathds{1}^TBRB^T\frac{\partial H}{\partial x}(x) +
\mathds{1}^TE\bar{d}=0,
\end{equation}
implying (as in the case $d=0$) that $\mathds{1}^Tx$ is a {\it conserved
quantity} for the closed-loop
system (\ref{closedloop}). In particular it
follows that the limit value $\lim_{t \to \infty}x(t) \in \spa\{ \mathds{1}\}$
is
determined by the initial condition $x(0)$.

\section{Basic setting with constrained flows}

In many cases of interest, the elements of the vector of flow inputs $u \in
\mathbb{R}^m$
corresponding to the edges of the graph will be {\it constrained}, that is
\begin{equation}
 u \in\mathcal{U}:=\{u\in\mathbb{R}^m\mid u^-\leqslant u\leqslant u^+\}
\end{equation}
for certain vectors $u^-$ and $ u^+$ satisfying $u^-_i<
u^+_i, i=1,\ldots,m$. In our previous \cite{Wei2013} we focused on
the case $u^-_i\leqslant0<u^+_i,i=1,2,\ldots,m.$ In the present paper we consider {\it arbitrary}
constraint intervals, necessitating a novel approach to the problem.

Thus we consider a general constrained version of the PI controller (\ref{PI}) 
discussed in the previous section, given as
\begin{equation}\label{PIconstrained}
\begin{array}{rcl}
\dot{x}_c & = & y ,\\[2mm]
u & = &\sat\big(-Ry - \frac{\partial H_c}{\partial x_c}(x_c)\,;u^-,u^+\big)
\end{array}
\end{equation}
For simplicity of exposition we consider throughout the rest of this paper the
standard Hamiltonian $H_c(x_c) = \frac{1}{2} \| x_c \|^2$ for the constrained PI
controller and the identity gain matrix $R=I$, while we throughout assume that
the Hessian matrix of Hamiltonian $H(x)$ is positive definite for any $x$.
Then the system (\ref{system1}) with nonzero in/outflows is given as
\begin{equation}\label{closedloop-sat-disturb}
\begin{aligned}
\dot{x} & =  B\sat\big(-B^T\frac{\partial H}{\partial
x}(x)-x_c\,;u^-,u^+\big)+E\bar{d},
\\[2mm]
\dot{x}_c & =  B^T\frac{\partial H}{\partial x}(x),
\end{aligned}
\end{equation}
In the rest of this section, we will show how the disturbance can be {\it
absorbed} into the constraint intervals and how the orientation can be made
compatible with the flow constraints.

First we note that we can incorporate the constant vector $\bar{d}$ of
in/outflows 
into the constraint intervals. Indeed, for
any $\eta\in\mathbb{R}^n$, we have the identity
\begin{equation}\label{identity}
\sat(x-\eta\,;u^-,u^+)+\eta=\sat(x\,;u^-+\eta,u^++\eta).
\end{equation}
Therefore for an in/out flow $\bar{d}$ satisfying the matching
condition, i.e., such that there exists $\bar{x}_c$ with $B\bar{x}_c=E\bar{d}$, we can
rewrite system (\ref{closedloop-sat-disturb}) as
\begin{equation}\label{disturbance_in_const}
 \begin{aligned}
  \dot{x} & = B\sat(-B^T\frac{\partial H}{\partial
x}(x)-\tilde{x}_c\,;u^-+\bar{x}_c,u^++\bar{x}_c), \\
\dot{\tilde{x}}_c & = B^T\frac{\partial H}{\partial
x}(x),
 \end{aligned}
\end{equation}
where $\tilde{x}_c=x_c-\bar{x}_c$.
It follows that, without loss of generality, we can restrict ourselves to the study of the closed-loop system
\begin{equation}\label{closedloop-sat}
\begin{aligned}
\dot{x} & =  B\sat\big(-B^T\frac{\partial H}{\partial
x}(x)-x_c\,;u^-,u^+\big),
\\[2mm]
\dot{x}_c & =  B^T\frac{\partial H}{\partial x}(x).
\end{aligned}
\end{equation}
for general $u^-$ and $ u^+$ with $u^-_i<
u^+_i, i=1,\ldots,m$ (where the vector $\bar{d}$ of in/outflows has been
incorporated in the vectors $u^-,u^+$).

\medskip

An essential ingredient in the analysis of the dynamical system
(\ref{closedloop-sat}) will be the following property of the scalar saturation
function $\sat(x;u^-,u^+)$, which allows us to split any edge in $\mathcal{G}$
into multiple edges. The scalar saturation function satisfies
\begin{equation}\label{division}
\begin{aligned}
 &\sat(x;u^-,u^+)\\
=&\sat(x;u^-,b_2)+\sum_{i=3}^{n-1}\sat(x-b_{i-1};0,b_i-b_{i-1})\\
&+\sat(x-b_{n-1};0,u^+-b_{n-1})
\end{aligned}
\end{equation}
for arbitrary $b_i,i=2,\ldots,n-1$ satisfying $u^-<b_2<\cdots<b_{n-1}<u^+.$
The above identity will imply that we can split any edge in $\mathcal{G}$ into 
multiple edges with the same orientation as the original one, and
with constraint
intervals $[u^-,b_2],[0,b_3-b_{2}],\ldots,[0,b_{n-1}-b_{n-2}],[0,u^+-b_{n-1}]$.
For any $i$-th edge in $\mathcal{G}$ the multiple edges resulting from
splitting 
of the $i$-th edge will be denoted as the $i_1$-th,\ldots,$i_{n-1}$-th edges.
Furthermore, we will denote the augmented graph which is generated by
splitting the $i$-th edge in $\mathcal{G}$ into multiple edges
by $\tilde{\mathcal{G}}$.

By choosing suitable initial conditions for the edge states at the newly added
edges of $\tilde{\mathcal{G}}$, the evolution of $x$ will be the same as that
in 
the original dynamical system (\ref{closedloop-sat}) defined on $\mathcal{G}$.
Indeed, corresponding to the identity
(\ref{division}) we can choose the initial conditions for the newly added edges as
follows
\begin{equation}
 \begin{aligned}
  x_{{c}_{i1}}(0)&=x_{ci}(0)\\
x_{{c}_{ik}}(0)&=x_{ci}(0)+b_k, k=2,\ldots,n-1,
 \end{aligned}
\end{equation}
where $x_{ci}(0)$ is the initial condition of the $i$-th edge state in
the dynamical
system (\ref{closedloop-sat}) defined on $\mathcal{G}$.

As a special case of the above construction, the {\it bi-directional} edge whose constraint interval satisfies
$u^-_i<0<u^+_i$, can be divided into {\it two uni-directional} edges with
constraint intervals $[u^-_i,0], [0,u^+_i]$ respectively, and the same
orientation.

Finally, we may change the {\it orientation} of
some of the
edges of the graph at will; replacing the corresponding columns $b_i$ of the
incidence matrix $B$ by $-b_i$. Noting the identity
\begin{equation}\label{identity1}
\sat(-x\,;u_i^-,u_i^+)=-\sat(x\,;-u_i^+,-u_i^-)
\end{equation}
this implies that we may
assume {\it without loss of generality} that the orientation of the graph is chosen
such that
\begin{equation}\label{assumption1}
u^+_i>0,\, i=1,2,\ldots,m.
\end{equation}
\begin{exmp}
Consider the graph given as in Fig.\ref{figure_ex1}, where the constraint
interval for edge $e_1$ is $[-2,-1]$. Clearly this network
is equivalent to the network where the edge direction is reversed from $v_2$ to
$v_1$ 
while the constraint interval is modified into $[1,2]$.
\end{exmp}
By dividing bi-directional edges into uni-directional ones and changing
orientations afterwards, we can also without loss of generality assume that
\begin{equation}\label{assumption2}
 u^-_i\geqslant0,\, i=1,2,\ldots,m.
\end{equation}
Conditions (\ref{assumption1}) and (\ref{assumption2}) will be standing assumptions throughout
the rest of the paper.
In general, we will say
that the graph is {\it compatible with the flow constraints} if (\ref{assumption1})
and (\ref{assumption2}) hold.

\section{Convergence conditions for the closed-loop dynamics with general flow constraints}

In this section, we will analyze system (\ref{closedloop-sat}) defined on a
general graph $\mathcal{G}$ with arbitrary constraint intervals. The main
construction is based on the following result which is proved in
\cite{Wei2013}.
\begin{lemma}\label{lemma}
A strongly connected graph is balanced if and only if it can
be covered by non-overlapping cycles.
\end{lemma}
The main idea for the subsequent analysis is now as follows. In view of Lemma 
\ref{lemma} the analysis of the system (\ref{closedloop-sat}) on a {\it
balanced} graph can be conducted
separately on each cycle. In other words, the behavior of the system
(\ref{closedloop-sat}) on a balanced graph is determined by the {\it subsystem}
defined on each cycle. Furthermore, these {\it subsystems} are independent of
each 
other, and as will follow from the subsequent Lemma \ref{cycle}, the steady
states of
the system (\ref{closedloop-sat}) defined on each cycle are determined only by the
constraint intervals. On the other hand, for a graph $\mathcal{G}$ which is {\it
not}
 balanced, we can split the overlapped
edges into multiple ones, using the construction explained in the previous
section, 
in order to {\it render} the graph balanced, and then use the same process as in
the balanced case.

Before delving into the analysis, let us consider two examples which show that
the stability of the system (\ref{closedloop-sat})
is dependent on the strong connectedness and on the constraint intervals,
especially the interval of the form $[0,u^+_i].$

\begin{exmp}\label{example1}
Consider the dynamical system (\ref{closedloop-sat}) defined on the
graph given by Fig.\ref{figure_ex1}
\begin{equation}
\begin{aligned}
\begin{bmatrix} \dot{x}_1 \\[2mm] \dot{x}_2 \end{bmatrix} &=
\begin{bmatrix} -1 \\[2mm] 1 \end{bmatrix}\sat(x_1-x_2-x_c,0,1)\\
\dot{x}_c &= x_2-x_1.
\end{aligned}
\end{equation}
This system will converge to a state satisfying $x_2>x_1$ and
$\sat(x_1-x_2-x_c,0,1)=0$. 
We see that although the graph $\mathcal{G}$ is not strongly connected, the
system still may reach a steady state.

\begin{figure}[ht]
\begin{center}
\begin{tikzpicture}
\tikzstyle{EdgeStyle}    = [thin,double= black,
                            double distance = 0.5pt]
\useasboundingbox (0,0) rectangle (4cm,0.5cm);
\tikzstyle{VertexStyle} = [shading         = ball,
                           ball color      = white!100!white,
                           minimum size = 20pt,%
                           inner sep       = 1pt,]
\Vertex[style={minimum
size=0.2cm,shape=circle},LabelOut=false,L=\hbox{$1$},x=1cm,y=0.2cm]{v1}
\Vertex[style={minimum
size=0.2cm,shape=circle},LabelOut=false,L=\hbox{$2$},x=4cm,y=0.2cm]{v2}
\draw
(v1) edge[->,>=angle 90,thin,double= black,double distance = 0.5pt]
node[above]{$e_1$} (v2);
\end{tikzpicture}
\caption{Illustrative graph}\label{figure_ex1}
\end{center}
\end{figure}
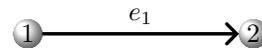
\end{exmp}

\begin{exmp}
Consider the dynamical system (\ref{closedloop-sat})
defined on the same graph as in Fig. \ref{figure_ex1}, but now with a different constraint
interval. The system can be written as
\begin{equation}
\begin{aligned}
\begin{bmatrix} \dot{x}_1 \\[2mm] \dot{x}_2 \end{bmatrix} &=
\begin{bmatrix} -1 \\[2mm] 1 \end{bmatrix}\sat(x_1-x_2-x_c,1,2)\\
\dot{x}_c &= x_2-x_1.
\end{aligned}
\end{equation}
At each time $t$, there will be positive flow from $x_1$ to $x_2$. Therefore
the states of system will go to plus or minus infinity. In this case, we call
the
system {\it unstable.}
\end{exmp}

As indicated in the beginning of this section, the analysis of the closed-loop system
(\ref{closedloop-sat}) defined on a cycle constitutes the cornerstone of the
analysis. The stability analysis on a cycle is given in the following lemma.

\begin{lemma}\label{cycle}
 Consider the closed-loop dynamical system  (\ref{closedloop-sat}) on a cycle whose orientation
is compatible with the constraint intervals $[u^-,u^+]$. The trajectories
of the closed-loop system (\ref{closedloop-sat}) converge to the set
\begin{equation}\label{set-zerodist}
 \mathcal{E}_{\mathrm{tot}} = \{ (x,x_c) \mid \frac{\partial H}{\partial x}(x) =
\alpha \mathds{1}_n, \, B\sat(-x_c\,;u^-,u^+) = 0 \}.
\end{equation}
if and only if the cycle is strongly connected and the intersection of all
the constraint intervals is again an interval with non-empty interior.
\end{lemma}

\begin{remark}
 Notice that when the graph contains cycles, the choice of $\bar{x}_c$ in
(\ref{disturbance_in_const}) is not
unique because for a cycle $\ker B = \spa\{\mathds{1}\}$. However, this fact does
not affect the condition in Lemma \ref{cycle}. Indeed, consider a cycle,
denoted as $\mathcal{C}$, whose
orientation is compatible with $[u^-,u^+]$ is strongly connected and
such that $\cap^{m}_{i=1}[u^-_i,u^+_i]$ has nonempty interior. Suppose that new
constraint intervals $[u^-+c\mathds{1},u^++c\mathds{1}]$ are imposed on
$\mathcal{C}$. If the orientation is compatible with the new constraint
intervals, then clearly $\mathcal{C}$ is strongly connected. If not, we can
prove that the cycle $\mathcal{C}'$ with reversed orientation with respect to
$\mathcal{C}$ is compatible with $[-u^+-c\mathds{1},-u^--c\mathds{1}]$ and again
strongly connected. Obviously, $\cap^{m}_{i=1}[u^-_i+c,u^+_i+c]$ and 
$\cap^{m}_{i=1}[-u^+_i-c,-u^-_i-c]$ both have nonempty
interiors. 
\end{remark}

\begin{proof}[of Lemma \ref{cycle}]
{\it Sufficiency}:
Consider the  Lyapunov
function given by
\begin{equation}\label{sat-Lyapunov}
V(x,x_c)=\mathds{1}^{T}_m S\big(-B^{T}\frac{\partial H}{\partial
x}(x)-x_c\,;u^-,u^+\big)+H(x),
\end{equation}
with
\begin{equation}
S(x\,;u^-,u^+)_i:=\int_0^{x_i} \sat(y\,;u^-_i,u^+_i)dy.
\end{equation}
The invariant set is given as
\begin{equation}
\begin{aligned}
 &\mathcal{I}=\{(\nu,x_c)\mid x_c=B^T\frac{\partial H}{\partial
x}(\nu)
t+x_c(0),\\
&B\sat\big(-B^T\frac{\partial H}{\partial
x}(\nu)-B^T\frac{\partial H}{\partial
x}(\nu)
t-x_c(0)\,;u^-,u^+\big)=0,\\
& \forall t\geq 0\}.
\end{aligned}
\end{equation}

For a strongly connected cycle, $\ker B = \spa\{\mathds{1}\}$. Suppose
$B^T\frac{\partial H}{\partial x}(\nu) \neq 0$, then there exists an edge, say
the $i$-th edge, whose flow reaches its upper bound, and an edge, say the
$j$-th 
edge whose flow reaches its lower bound.  Because $[u^-_i,u^+_i]$ and
$[u^-_j,u^+_j]$ are overlapped, it follows that
\begin{equation}
 u^+_i>u^-_j
\end{equation}
Then the vector whose $i$-th component is $u^+_i$ and $j$-th component is
$u^-_j$ 
does not belong to $\spa\{\mathds{1}\}$. Therefore, for $t$ large enough,
\begin{equation}
 B\sat\big(-B^T\frac{\partial H}{\partial
x}(\nu)-B^T\frac{\partial H}{\partial
x}(\nu)
t-x_c(0)\,;u^-,u^+\big)\neq 0
\end{equation}
and we have reached a contradiction.

{\it Necessity}:
First, suppose that the cycle compatible with the constraint interval is not strongly
connected. Say there is a path from $x_i$ to
$x_j$, but not a path from $x_j$ to $x_i$. In other words, there can be
a positive flow from $x_i$ to $x_j$, but not vice versa. Then for suitable
initial conditions, $\frac{\partial H}{\partial x_i} (x(t)) <\frac{\partial H}{\partial
x_j}(x(t))$ for all $t\geqslant0.$

Secondly, suppose the graph compatible with constraints interval is strongly
connected, but there exist two constraints intervals such that their
intersection is empty, then the system (\ref{closedloop-sat}) is unstable.
Indeed, suppose $[u^-_i, u^+_i]\cup[u^-_j,u^+_j]=\emptyset$, where, without loss
of generality, we can assume $u^-_i>u^+_j$. So there will be more positive flow
along the $i$-th edge than along the $j$-th edge, which makes the system unstable.

Now we analyze the case that the intersection of any two constraints intervals
is not empty but a single point. Without loss of generality,
\begin{equation}
 [u^-_i,u^+_i]\cap[u^-_j,u^+_j]=\{u^+_i\}
\end{equation}
and $u^+_i\in[u^-_k,u^+_k], k=1,2,\cdots,m.$ So there exist
$B^T\frac{\partial H}{\partial x}(\nu)\neq0$ and suitable $x_c(0)$ such that
\begin{equation}
 B\sat\big(-B^T\frac{\partial H}{\partial
x}(\nu)-B^T\frac{\partial H}{\partial x}(\nu)t-x_c(0)\,;u^-,u^+\big)=0,
\end{equation}
for all $t\geqslant0$, that is
\begin{equation}
 \sat\big(-B^T\frac{\partial H}{\partial
x}(\nu)-B^T\frac{\partial H}{\partial
x}(\nu)t-x_c(0)\,;u^-,u^+\big)=u^+_i\mathds{1}.
\end{equation}
In this case, $\nu$ is an equilibrium for $x$ satisfying
$B^T\frac{\partial H}{\partial x}(\nu)\neq0$. In fact, flows in those edges
which belong to $\mathcal{E}_1=\{k \mbox{-th edge} \mid u^+_k=u^+_i\}$ reach
their upper bounds, while flows in the edges
which belong to $\mathcal{E}_2=\{k \mbox{-th edge} \mid u^-_k=u^+_i\}$ reach
their lower bounds. Thus $\frac{\partial H}{\partial x}$ will form a clustering,
and no consensus will be reached
\end{proof}

\begin{corollary}
The state $x$ will converge to a clustering if and only if the intersection of all the
constraint intervals is only a single point. The system is unstable if the
intersection of all the constraint intervals is empty.
\end{corollary}

\begin{exmp}\label{example}
Consider the dynamical system (\ref{closedloop-sat}) defined on the
Fig.\ref{figure}. We will show three different constraints intervals and the
corresponding results.

1. The constraint intervals for the edges $e_1,e_2,e_3$ are $[1,2], [2,3],
[0,3]$ 
respectively. In this case $x$ will converge to a clustering. The result is
given in
Fig.\ref{clustering}.

2. If we consider constraint intervals $[1,2.5], [2,3], [0,3]$ for the edges
$e_1,e_2,e_3$, 
then $x$ will converge to consensus, as can be seen from
Fig.\ref{consensus}.

3. Suppose the constraint intervals for $e_1,e_2,e_3$ are $[1,1.5], [2,3], [0,3]$
respectively. In this case $x$ will explode. The result is given in
Fig.\ref{unstable}.
\begin{figure}[ht]
\begin{center}
\begin{tikzpicture}
\tikzstyle{EdgeStyle}    = [thin,double= black,
                            double distance = 0.5pt]
\useasboundingbox (0,0) rectangle (3cm,1.5cm);
\tikzstyle{VertexStyle} = [shading         = ball,
                           ball color      = white!100!white,
                           minimum size = 20pt,%
                           inner sep       = 1pt,]
\Vertex[style={minimum
size=0.2cm,shape=circle},LabelOut=false,L=\hbox{$1$},x=1.5cm,y=1.5cm]{v1}
\Vertex[style={minimum
size=0.2cm,shape=circle},LabelOut=false,L=\hbox{$2$},x=0.5cm,y=0cm]{v2}
\Vertex[style={minimum
size=0.2cm,shape=circle},LabelOut=false,L=\hbox{$3$},x=2.5cm,y=0cm]{v3}
\draw
(v1) edge[->,>=angle 90,thin,double= black,double distance = 0.5pt]
node[left]{$e_1$} (v2)
(v2) edge[->,>=angle 90,thin,double= black,double distance = 0.5pt]
node[above]{$e_2$} (v3)
(v3) edge[->,>=angle 90,thin,double= black,double distance = 0.5pt]
node[right]{$e_3$} (v1);
\end{tikzpicture}
\caption{Network of Example $\ref{example}$}\label{figure}
\end{center}
\end{figure}
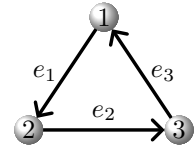

\begin{figure}[htbp]
\centering
\subfigure[Clustering]{
 \includegraphics[height=2.9cm]{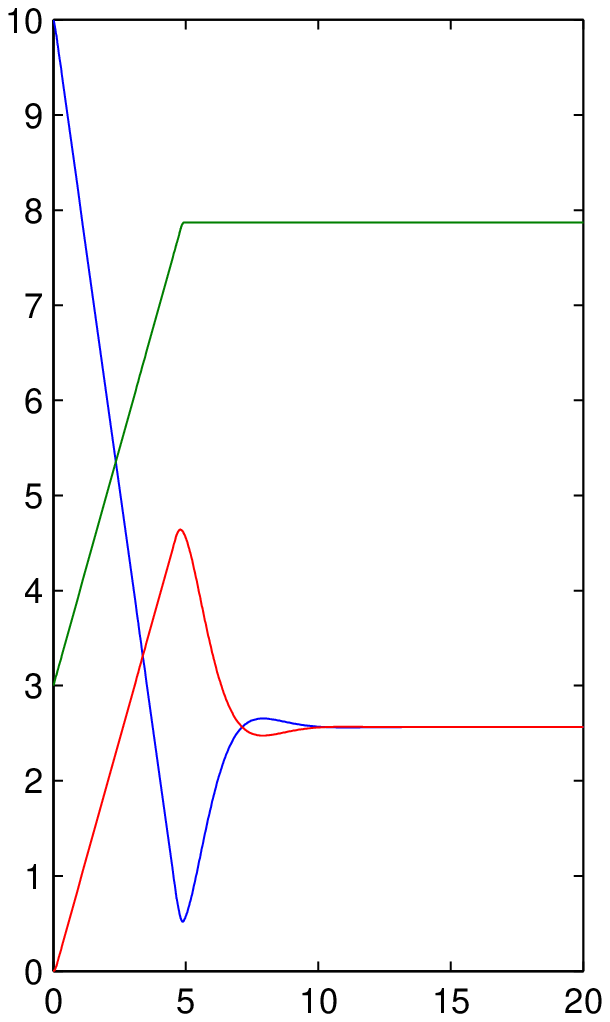}
 \label{clustering}}
\subfigure[Consensus]{
 \includegraphics[height=2.9cm]{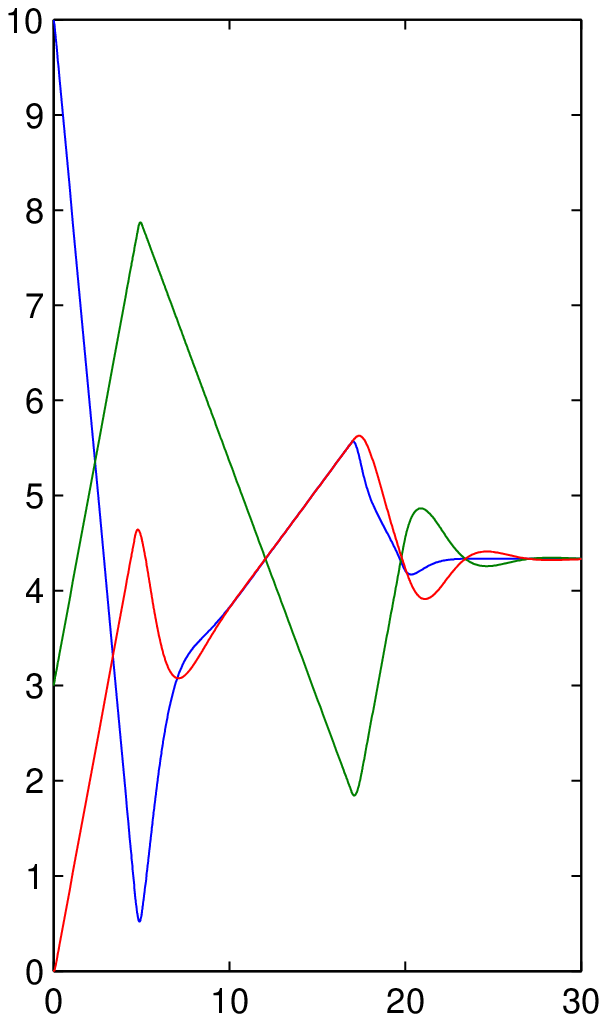}
 \label{consensus}}
\subfigure[Unstable]{
 \includegraphics[height=2.9cm]{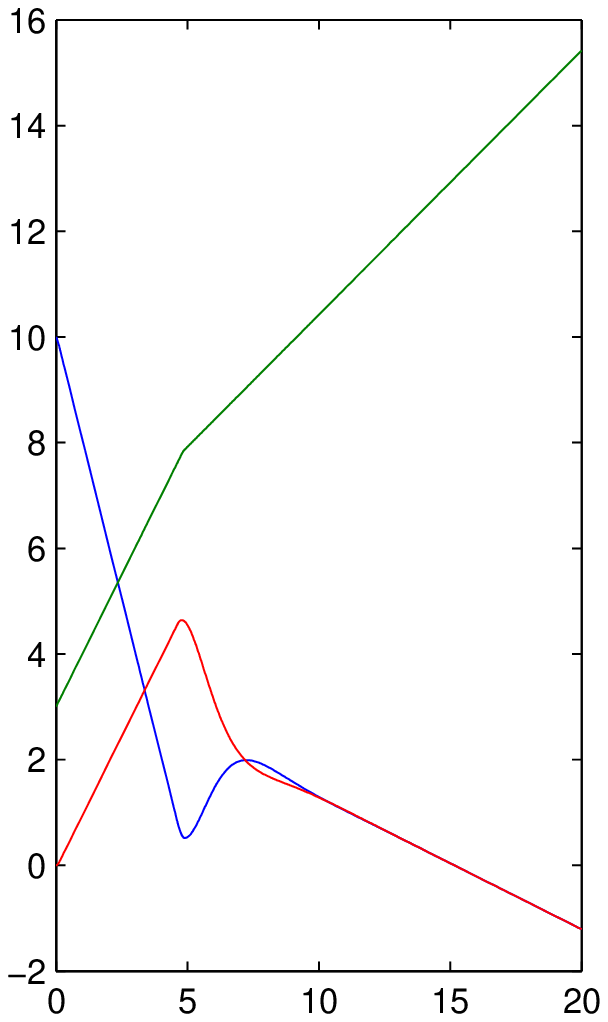}
 \label{unstable}}
\caption{The trajectories of the storage at the vertices}\label{example3}
 \end{figure}
\end{exmp}

Now let us consider the closed-loop system (\ref{closedloop-sat}) defined on
a general graph. At this moment we will only give a {\it sufficient condition} for the system
(\ref{closedloop-sat}) under arbitrary constraints to reach load balancing
(consensus). Consider a strongly connected network compatible with
the constraint intervals $[u^-,u^+]$. According to Lemma \ref{lemma}, suppose
there exists $k$ cycles to cover the graph, denoted as 
$\mathcal{T}=(C_1,C_2,\ldots,C_k)$. Given
$\mathcal{T}$, we can define a {\it multiplicity vector} $T\in\mathbb{R}^m$
whose $i$-th component is the number of cycles in $\mathcal{T}$ which contain the
$i$-th edge. Then we construct an augmented network
$\tilde{\mathcal{G}}(\mathcal{T})$ by splitting each edge of $\mathcal{G}$ into
multiple edges based on their multiplicities, using the identity (\ref{division}).
For instance, if the $i$-th edge of $\mathcal{G}$ has been used $T_i$ times in
$\mathcal{T}$ then we splitting $i$-th edge into $T_i$ edges. The newly generated
edges have constraint intervals
$[u^-_i,b_2],[0,b_3-b_2],\ldots,[0,u^+_i-b_{T_i}],$ for arbitrary
$u^-_i<b_2<\cdots<b_{T_i}<u^+_i$. Furthermore, it can be easily seen that
$\tilde{\mathcal{G}}(\mathcal{T})$ is balanced, and that it can be covered by
non-overlapping cycles.  We denote the set of cycles to cover
$\tilde{\mathcal{G}}(\mathcal{T})$ by $\tilde{\mathcal{T}}$. The above process
can be explained by the following example.

\begin{exmp}
In this example, we consider the graph $\mathcal{G}$ given as in
Figure. \ref{combined}(left). Notice that $\mathcal{G}$ is unbalanced and
that $\mathcal{T}=\{C_1,C_2\}$ is a minimal covering set where
$C_1=\{e_1,e_2,e_3\}$ and $C_1=\{e_3,e_4,e_5\}$. So the corresponding
multiplicity vector $T$ is given as $T=[1,1,2,1,1]^T.$

By dividing $e_3$ into two edges, we obtain the augmented graph
$\tilde{\mathcal{G}}(\mathcal{T})$ as in Figure.\ref{combined}(right). Here the
constraint
intervals for the edge $e_{3_1}$ and the edge $e_{3_2}$ in
$\tilde{\mathcal{G}}(\mathcal{T})$ are $[u^-_3,b],[0,u^+_3-b]$ respectively,
while $[u^-_3,u^+_3]$ is the constraint interval for $e_3$ in $\mathcal{G}$ with $u^-_3<b<u^+_3$.

Now $\tilde{\mathcal{G}}(\mathcal{T})$ is balanced and can be covered by
non-overlapping cycles. Indeed,
$\tilde{\mathcal{T}}=\{\tilde{C}_1,\tilde{C}_2\}$ where
$\tilde{C}_1=\{e_1,e_2,e_{3_1}\}$ and $\tilde{C}_2=\{e_4,e_5,e_{3_2}\}.$

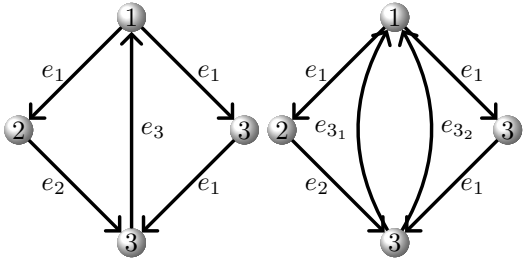
\begin{figure}
\centering
\begin{tikzpicture}
\tikzstyle{EdgeStyle}    = [thin,double= black,
                            double distance = 0.5pt]
\useasboundingbox (0,0) rectangle (7cm,3cm);
\tikzstyle{VertexStyle} = [shading         = ball,
                           ball color      = white!100!white,
                           minimum size = 20pt,%
                           inner sep       = 1pt,]
\Vertex[style={minimum
size=0.2cm,shape=circle},LabelOut=false,L=\hbox{$1$},x=1.5cm,y=3cm]{v1}
\Vertex[style={minimum
size=0.2cm,shape=circle},LabelOut=false,L=\hbox{$2$},x=0cm,y=1.5cm]{v2}
\Vertex[style={minimum
size=0.2cm,shape=circle},LabelOut=false,L=\hbox{$3$},x=1.5cm,y=0cm]{v3}
\Vertex[style={minimum
size=0.2cm,shape=circle},LabelOut=false,L=\hbox{$3$},x=3cm,y=1.5cm]{v4}
\Vertex[style={minimum
size=0.2cm,shape=circle},LabelOut=false,L=\hbox{$1$},x=5cm,y=3cm]{v5}
\Vertex[style={minimum
size=0.2cm,shape=circle},LabelOut=false,L=\hbox{$2$},x=3.5cm,y=1.5cm]{v6}
\Vertex[style={minimum
size=0.2cm,shape=circle},LabelOut=false,L=\hbox{$3$},x=5cm,y=0cm]{v7}
\Vertex[style={minimum
size=0.2cm,shape=circle},LabelOut=false,L=\hbox{$3$},x=6.5cm,y=1.5cm]{v8}
\draw
(v1) edge[->,>=angle 90,thin,double= black,double distance = 0.5pt]
node[left]{$e_1$} (v2)
(v2) edge[->,>=angle 90,thin,double= black,double distance = 0.5pt]
node[left]{$e_2$} (v3)
(v3) edge[->,>=angle 90,thin,double= black,double distance = 0.5pt]
node[right]{$e_3$} (v1)
(v1) edge[->,>=angle 90,thin,double= black,double distance = 0.5pt]
node[right]{$e_1$} (v4)
(v4) edge[->,>=angle 90,thin,double= black,double distance = 0.5pt]
node[right]{$e_1$} (v3)
(v5) edge[->,>=angle 90,thin,double= black,double distance = 0.5pt]
node[left]{$e_1$} (v6)
(v6) edge[->,>=angle 90,thin,double= black,double distance = 0.5pt]
node[left]{$e_2$} (v7)
(v7) edge[->,>=angle 90,bend left,thin,double= black,double distance = 0.5pt]
node[left]{$e_{3_1}$} (v5)
(v7) edge[->,>=angle 90,bend right,thin,double= black,double distance = 0.5pt]
node[right]{$e_{3_2}$} (v5)
(v5) edge[->,>=angle 90,thin,double= black,double distance = 0.5pt]
node[right]{$e_1$} (v8)
(v8) edge[->,>=angle 90,thin,double= black,double distance = 0.5pt]
node[right]{$e_1$} (v7);
\end{tikzpicture}

\caption{Left: The graph $\mathcal{G}$. Right: The augmented graph
$\tilde{\mathcal{G}}(\mathcal{T})$.The generation of the augmented graph
$\tilde{\mathcal{G}}(\mathcal{T})$ based on $\mathcal{T}$.}
\label{combined}
\end{figure}

\end{exmp}

The main result of the paper can be summarized as the following theorem substantially generalizing Lemma \ref{cycle}.

\begin{theorem}\label{main}
Consider the closed-loop dynamical system (\ref{closedloop-sat}) defined on a strongly connected
graph which is compatible with the constraint intervals. Let
$\mathcal{T}$ be a
minimal covering set for $\mathcal{G}$ and let $\tilde{\mathcal{G}}(\mathcal{T})$
be the augmented graph based on $\mathcal{T}$. Let
$\tilde{\mathcal{T}}=(\tilde{C}_1,\tilde{C}_2,\ldots,\tilde{C}_k)$ be a covering
set of cycles for $\tilde{\mathcal{G}}(\mathcal{T})$. If there exists a splitting
of the overlapped edges in $\mathcal{G}$ such that the intersection of all constraint
intervals of each cycle $\tilde{C}_i,i=1,2,\ldots,k$ has non-empty interior, then the trajectories of the system
(\ref{closedloop-sat}) will converge to
\begin{equation}
\begin{aligned}
 \mathcal{E}_{\mathrm{tot}} = & \{ (x,x_c) \mid \frac{\partial H}{\partial x}(x) =
\alpha \mathds{1}, \, \alpha \in \mathbb{R}, \,\\
& B\sat(-x_c\,;u^-,u^+)=0
\, \}.
\end{aligned}
\end{equation}
\end{theorem}
\medskip
\begin{proof}
Because of lack of space, we only give a sketch of the proof. Consider the
same Lyapunov function (\ref{sat-Lyapunov}). If we choose a constant vector
$(\nu,x_c(0))\in\mathcal{I}$, which is the largest invariant set in
$\{(x,x_c)\mid \dot{V}=0\}$, then along this trajectory $V(\nu,
B^T\frac{\partial H}{\partial x}(\nu)t+x_c(0))$ is constant for all time
$t\geqslant0$. 
Suppose $B^T\frac{\partial H}{\partial x}(\nu)\neq 0$, then by the fact that
$\tilde{\mathcal{G}}(\mathcal{T})$ can be covered by non-overlapping cycles, we
can prove that for $t$ large enough, $\frac{d}{dt}V(\nu,
B^T\frac{\partial H}{\partial x}(\nu)t+x_c(0))>0$. This yields a contradiction.
\end{proof}

\begin{exmp}
The sufficiency condition in Theorem \ref{main} is not a necessary condition. 
Indeed, consider the dynamic (\ref{closedloop-sat}) defined on the
network given in Fig.\ref{combined}(left). The constraint intervals for
$e_i,i=1,2,\ldots,5$ are $[0.3,1]$, $[0.3,1]$, $[0.5,0.8]$, $[0.3,1]$, $[0.3,1]$
respectively. There does not exist any splitting such that the intersections of the constraint intervals have
nonempty interior. However $\frac{\partial H}{\partial x}(x(t))$ converges to
consensus. A special case with $H(x)=\frac{1}{2}\|x\|^2$ is shown in
Fig.\ref{final_example}.
\begin{figure}[htbp]
\centering
 \includegraphics[height=2cm,width=6.8cm]{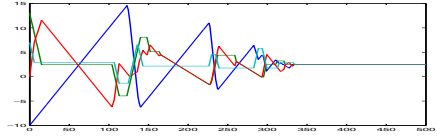}
\caption{The trajectories of the storage at the vertices}\label{final_example}
 \end{figure}
\end{exmp}

\section{CONCLUSIONS}

We have discussed a basic model of dynamical distribution networks where the
flows through the edges are generated by distributed PI controllers. The
resulting system can be naturally modeled as a port-Hamiltonian system with
arbitrary 
flow constraint intervals. Key tools
in the analysis are the construction of a $C^1$ Lyapunov function and the
observation given in Lemma \ref{lemma}. Based on that, we have derived
necessary and sufficient conditions for asymptotic consensus and clustering for
a dynamical system defined on a cycle. For arbitrary networks we have obtained a sufficient
condition for consensus or clustering.

An obvious open problem is to find sufficient and necessary conditions for
an arbitrary network to reach consensus or clustering. This is currently under
investigation.
Many other questions can be addressed within the same framework. For example, what is happening if
the in/outflows are not assumed to be constant, but are e.g. periodic
functions of time; see already \cite{depersis}.

\addtolength{\textheight}{-12cm}   





\bibliographystyle{IEEEtran} 
\bibliography{ifacconf}

\begin{thebibliography}{10}
\providecommand{\url}[1]{#1}
\csname url@rmstyle\endcsname
\providecommand{\newblock}{\relax}
\providecommand{\bibinfo}[2]{#2}
\providecommand\BIBentrySTDinterwordspacing{\spaceskip=0pt\relax}
\providecommand\BIBentryALTinterwordstretchfactor{4}
\providecommand\BIBentryALTinterwordspacing{\spaceskip=\fontdimen2\font plus
\BIBentryALTinterwordstretchfactor\fontdimen3\font minus
  \fontdimen4\font\relax}
\providecommand\BIBforeignlanguage[2]{{%
\expandafter\ifx\csname l@#1\endcsname\relax
\typeout{** WARNING: IEEEtran.bst: No hyphenation pattern has been}%
\typeout{** loaded for the language `#1'. Using the pattern for}%
\typeout{** the default language instead.}%
\else
\language=\csname l@#1\endcsname
\fi
#2}}

\bibitem{Wei2012}
{A.J. van der Schaft and J. Wei}, ``A hamiltonian perspective on the control of
  dynamical distribution networks,'' \emph{4th IFAC Workshop on Lagrangian and
  Hamiltonian Methods for Non Linear Control}, pp. 24--29, 2012.

\bibitem{schaftSIAM}
{A.J. van~der Schaft and B.M. Maschke}, ``Port-{H}amiltonian systems on
  graphs,'' \emph{SIAM J. Control and Optimization}, vol. 51(2), pp. 906--937,
  2013.

\bibitem{schaftCDC08}
{A.J. van der Schaft and B.M. Maschke}, ``Conservation laws on
  higher-dimensional networks,'' \emph{Proc. 47th IEEE Conf. on Decision and
  Control}, 2008.

\bibitem{schaftBosgrabook}
------, \emph{Model-Based Control: Bridging Rigorous Theory and Advanced
  Technology, P.M.J. Van den Hof, C. Scherer, P.S.C. Heuberger, eds., chapter
  Conservation laws and lumped system dynamics}.\hskip 1em plus 0.5em minus
  0.4em\relax Berlin-Heidelberg: Springer, 2009.

\bibitem{schaftNECSYS10}
------, ``Port-{H}amiltonian dynamics on graphs: Consensus and coordination
  control algorithms,'' \emph{Proc. 2nd IFAC Workshop on Distributed Estimation
  and Control in Networked Systems}, pp. 175--178, 2010.

\bibitem{allgower11}
{M. B\"{u}rger and D. Zelazo and F. Allg\"{o}wer}, ``Network clustering: A
  dynamical systems and saddle-point perspective,'' \emph{IEEE Conference on
  Decision and Control}, pp. 7825--7830, 2011.

\bibitem{Mesbahi11}
{D. Zelazo and M. Mesbahi}, ``Edge agreement: Graph-theoretic performance
  bounds and passivity analysis,'' \emph{Automatic Control, IEEE Transactions
  on}, vol.~56, no.~3, pp. 544 --555, march 2011.

\bibitem{Blanchini2000}
F.~Blanchini, S.Miani, and W.Ukovich, ``Control of production-distribution
  systems with unknown inputs and system failures,'' \emph{Automatic Control,
  IEEE Transactions on}, vol.~45, no.~6, pp. 1072--1081, 2000.

\bibitem{Bauso2011}
D.Bauso, F.Blanchini, L.Giarr\'{e}, and R.Pesenti., ``A decentralized solution
  for the constrained minimum-norm flow,'' \emph{Submitted to Automatic
  Control, IEEE Transactions on}, 2011.

\bibitem{Wei2013}
{J. Wei and A.J. van der Schaft}, ``Load balancing of dynamical distribution
  networks with flow constraints and unknown in/outflows,'' \emph{Accepted by
  Systems \& Control Letters}, 2013.

\bibitem{Bollobas98}
{B. Bollobas}, \emph{Modern Graph Theory}, ser. Graduate Texts in
  Mathematics.\hskip 1em plus 0.5em minus 0.4em\relax New York: Springer, 1998,
  vol. 184.

\bibitem{vanderschaftmaschkearchive}
{A.J van~der Schaft and B.M. Maschke}, ``The {H}amiltonian formulation of
  energy conserving physical systems with external ports,'' \emph{Archiv
  f\"{u}r Elektronik und \"{U}bertragungstechnik}, vol.~49, pp. 362--371, 1995.

\bibitem{vanderschaftbook}
{A.J. van der Schaft}, \emph{$L_2$-Gain and Passivity Techniques in Nonlinear
  Control}, ser. Lect. Notes in Control and Information Sciences.\hskip 1em
  plus 0.5em minus 0.4em\relax Berlin: Springer-Verlag, 1996, vol. 218.

\bibitem{depersis}
C.~D. Persis, ``Balancing time-varying demand-supply in distribution networks:
  an internal model approach,'' \emph{arxiv}, 2013.

\end{thebibliography}

\end{document}